\documentclass[psamsfonts]{amsart}

\usepackage{amssymb,amsfonts}
\usepackage[all,arc]{xy}
\usepackage{enumerate}
\usepackage{mathrsfs}
\usepackage{graphicx}

\newtheorem{theorem}{Theorem}[section]
\newtheorem{corollary}[theorem]{Corollary}
\newtheorem{proposition}[theorem]{Proposition}
\newtheorem{lemma}[theorem]{Lemma}

\theoremstyle{definition}
\newtheorem{definition}[theorem]{Definition}

\newtheorem{example}[theorem]{Example}

\newtheorem{remark}[theorem]{Remark}

\DeclareMathOperator{\Id}{Id}
\DeclareMathOperator{\Max}{Max}

\DeclareMathOperator{\gp}{gp}


\numberwithin{equation}{section}

\begin{document}

\title[Valuation Semirings]{Valuation Semirings}

\author[Peyman Nasehpour]{\bfseries Peyman Nasehpour}

\address{Peyman Nasehpour\\
Department of Engineering Science \\
Faculty of Engineering \\
University of Tehran   \\
Tehran\\
Iran}
\email{nasehpour@gmail.com}



\begin{abstract}
The main scope of this paper is to introduce valuation semirings in general and discrete valuation semirings in particular. In order to do that, first we define valuation maps and investigate them. Then we define valuation semirings with the help of valuation maps and prove that a multiplicatively cancellative semiring is a valuation semiring if and only if its ideals are totally ordered by inclusion. We also prove that if the unique maximal ideal of a valuation semiring is subtractive, then it is integrally closed. We end this paper by introducing discrete valuation semirings and show that a semiring is a discrete valuation semiring if and only if it is a multiplicatively cancellative principal ideal semiring possessing a nonzero unique maximal ideal. We also prove that a discrete valuation semiring is a Gaussian semiring if and only if its unique maximal ideal is subtractive.
\end{abstract}

\maketitle

\centerline{In memory of Prof. Dr. Manfred Kudlek}

\section{Introduction}

Semirings not only have significant applications in different fields such as automata theory in theoretical computer science, (combinatorial) optimization theory, and generalized fuzzy computation, but are fairly interesting generalizations of two broadly studied algebraic structures, i.e., rings and bounded distributive lattices. Valuation theory for rings was introduced by Krull in \cite{Kr} and has been proven to be a very useful tool in ring theory. Also valuation theory represents a nice interplay between ring theory and ordered Abelian groups. The main purpose of this paper is to generalize valuation theory for semirings.\\

In this paper, all semirings are commutative with a zero and a nonzero identity. In the first section of the present paper, similar to the classical concept of valuations in Bourbaki (\cite[VI, 3.1]{Bo}), we define the concept of valuation maps for semirings with values in totally ordered commutative monoids (tomonoids for short). Note that our notion for valuation maps should not to be confused with the concept of valuation functions for semirings in \cite[Definition 6.4]{HW}.

By \emph{an $M$-valuation on $S$}, we mean a map $v: S \rightarrow M_{\infty}$ with the following properties:

\begin{enumerate}

\item $S$ is a semiring and $M_{\infty}$ is a tomonoid with the greatest element $+\infty$, which has been obtained from the tomonoid $M$ with no greatest element,

\item $v(xy) = v(x) + v(y)$, for all $x,y\in S$,

\item $v(x+y) \geq \min\{v(x),v(y)\}$, for all $x,y\in S$,

\item $v(1) = 0$ and $v(0) = +\infty$.

\end{enumerate}

When there is an $M$-valuation $v$ on $S$, it is easy to see that the set $S_v = \{s\in S : v(s)\geq 0 \}$ is a subsemiring of $S$ and plays an important role in our paper. If the map $v$ is surjective, we show that the semiring $S$ is a semifield if and only if $M$ is an Abelian group and the set of units $U(S_v)$ of $S_v$ is equal to $\{s\in S_v : v(s)=0 \}$ (See Theorem \ref{units2}).

Let us recall that a nonempty subset $I$ of a semiring $S$ is said to be an ideal of $S$, if $a,b \in I$ implies $a+b \in I$ and $sa\in I$ for any $s\in S$. An ideal $I$ of a semiring $S$ is said to be subtractive (in some references k-ideal), if $a+b \in I$ and $a \in I$ implies $b \in I$ for all $a,b \in S$. At last, an ideal $I$ of a semiring $S$ is said to be proper, if $I \neq S$.

Let $v$ be an $M$-valuation on $S$. We define $v$ to have $\min$-\emph{property}, whenever $v(x) \neq v(y)$, then $v(x+y) = \min\{v(x),v(y)\}$ for any $x,y \in S$. Though this property holds for valuation maps on rings, this is not the case for semirings (Check Remark \ref{min-property}). We say ``$T$ is a $V$-semiring with respect to the triple $(S,v,M)$'', if $T$ is a semiring and there exists an $M$-valuation $v$ on $S$ such that the semiring $S$ contains $T$ as a subsemiring and $T=S_v = \{s\in S : v(s)\geq 0 \} $ (Definition \ref{Vsemiring}). Actually in Corollary \ref{subtractive2}, we show that if $T$ is a $V$-semiring with respect to the triple $(K,v,M)$, where $K$ is a semifield and $v$ is a surjective map, then the following statements are equivalent:

\begin{enumerate}

\item The valuation map $v$ has the $\min$-property;

\item The maximal $K_v$-ideal $P_v = \{x\in S_v : v(x) > 0 \}$ is subtractive.

\end{enumerate}

Similar to ring theory, an ideal $A$ of a semiring $S$ is said to be principal if $A$ is generated by one element, i.e., $A = (a) = \{sa : s\in S\}$ for some $a\in S$. In Proposition \ref{principal}, we prove that if $S_v$ is a $V$-semiring with respect to the triple $(K,v,M)$, where $K$ is a semifield and $v$ is surjective, then the principal ideals of $S_v$ are totally ordered by inclusion. This is the base for the definition of valuation semirings given in section 2.

Let us recall a notation. If $S$ is a multiplicatively cancellative semiring (for short MC semiring), we denote its semifield of fractions by $F(S)$ (Refer to Remark \ref{semifieldoffractions}). In section 2, we define a semiring $S$ to be a \emph{valuation semiring}, if there exists an $M$-valuation $v$ on $K$, where $K$ is a semifield containing $S$ as a subsemiring, $v$ is surjective, and $S = K_v = \{s\in K : v(s)\geq 0 \}$. One of the main theorems of the paper is the following (Theorem \ref{valuationsemiring}):

An MC semiring $S$ is a valuation semiring, if and only if one of the following equivalent conditions hold:

\begin{enumerate}

\item For any element $x\in F(S)$, either $x\in S$ or $x^{-1}\in S$,

\item For any ideals $I,J$ of $S$, either $I \subseteq J$ or $J \subseteq I$,

\item For any elements $x,y\in S$, either $(x) \subseteq (y)$ or $(y) \subseteq (x)$.

\end{enumerate}

Let $S$ be an MC semiring and $F(S)$ its semifield of fractions. The element $u\in F(S)$ is said to be integral over $S$ if there exist $a_1,a_2,\ldots, a_n$ and $b_1,b_2,\ldots,b_n$ in $S$ such that $u^n+a_1 u^{n-1}+\cdots + a_n = b_1 u^{n-1}+ \cdots + b_n$. The semiring $S$ is said to be integrally closed if the set of elements of $F(S)$ that are integral over $S$ is the set $S$ \cite[p. 88]{DM}.

In Theorem \ref{integralthm}, we show that if $S$ is a valuation semiring such that its unique maximal ideal $J$ is subtractive, then $S$ is integrally closed.

Section 3 is devoted to discrete valuation semirings. We define a semiring $S$ to be a \emph{discrete valuation semiring} (DVS for short), if $S$ is a $V$-semiring with respect to the triple $(K,v,\mathbb Z)$, where $K$ is a semifield and $v$ is surjective (Definition \ref{discrete}). Before explaining briefly what we do in section 3, we need to recall some concepts. A semiring is said to be quasi-local if it has only one maximal ideal. A semiring satisfies ACCP if any ascending chain of principal ideals of $S$ is stationary. A semiring $S$ is said to be a principal ideal semiring, if each ideal of $S$ is principal.

In this section, we characterize discrete valuation semirings as follows (See Theorem \ref{ACCP}): A semiring $S$ is a discrete valuation semiring if and only if one of the following equivalent conditions satisfies:

\begin{enumerate}

\item $S$ is a principal ideal MC semiring possessing a unique maximal ideal $J \neq (0)$,

\item $S$ is a quasi-local MC semiring whose unique maximal ideal $J \neq (0)$ is principal and $\bigcap_{n=1}^{\infty} J^n = (0)$,

\item $S$ is a quasi-local MC semiring whose unique maximal ideal $J \neq (0)$ is principal, which satisfies ACCP.

\end{enumerate}

Let us recall that if $S$ is a semiring, for a polynomial $f \in S[X]$, the content of $f$, denoted by $c(f)$, is defined to be the $S$-ideal generated by the coefficients of $f$. A semiring $S$ is called Gaussian if $c(fg)=c(f)c(g)$ for all polynomials $f,g \in S[X]$ (\cite[Definition 7]{N}). Finally in Theorem \ref{min-property2}, we prove the following:

Let $S$ be a discrete valuation semiring. Then the following statements are equivalent:

\begin{enumerate}

\item The unique maximal ideal of $S$ is subtractive,

\item The valuation map $v : F(S) \rightarrow \mathbb Z_{\infty}$ satisfies the $\min$-property,

\item Each ideal of the semiring $S$ is subtractive,

\item The semiring $S$ is Gaussian.

\end{enumerate}

A connoisseur of valuation theory in commutative algebra already have noticed that some of the definitions and results, mentioned above, are extensions of their ring version ones. Therefore, throughout the paper, some examples are given to show that some objects really satisfy the hypotheses of certain definitions and theorems, while the others are given to show that the hypotheses of certain theorems in the paper cannot be greatly weakened.

Though we have tried our paper to be self-contained in definitions and terminology, but many of them can be found in the book \cite{Go2}. One may refer to the books \cite{HW}, \cite{Go2} and \cite{P} for more on semirings.

\section{Valuation Maps on Semirings}

The main scope of this section is to generalize the concept of valuation maps on rings and investigate the properties of these maps. While according to the classical definition of valuation maps in Bourbaki (\cite[VI, 3.1]{Bo}), the values of the elements of a ring belong to a totally ordered Abelian group, we define valuation maps in such way that the values of the elements of a semiring belong to an arbitrary totally ordered commutative monoid. We recall that by a totally ordered commutative monoid (tomonoid for short) $(M,+,0,\leq)$, it is meant a commutative monoid $(M,+,0)$ such that $(M,\leq)$ is a totally ordered set and $x\leq y$ implies $x+z \leq y+z$ for any $z\in M$. We also need to recall the following:

In monoid theory, there is a routine technique that one can obtain a monoid with the greatest element from a monoid with no greatest element:

Let $(M,+,0,\leq)$ be a tomonoid with no greatest element. One can configure a set $M_{\infty}$ by adjoining a new element to the monoid $M$, denoted by $+\infty$, in this way that the total order on $M$ induces a total order on $M_{\infty}$ such that $+\infty$ is the greatest element, i.e., for all $m\in M$, we have that $m < +\infty$ and the monoid structure on $M$ induces a monoid structure on $M_{\infty}$ with the following rule: $$m + (+\infty) = (+\infty) + m = +\infty, \forall~ m\in M_{\infty}.$$

It is immediate that $M_{\infty}$ is also a tomonoid with the greatest element $+\infty$. Now similar to the classical concept of valuations in Bourbaki (\cite[VI, 3.1]{Bo}), we define the concept of valuation maps for semirings with values in tomonoids (not to be confused with the concept of valuation functions for semirings in \cite[Definition 6.4]{HW}).

\begin{definition}

\label{valuation}

By ``\emph{an $M$-valuation on $S$}", we mean a map $v: S \rightarrow M_{\infty}$ with the following properties:

\begin{enumerate}

\item $S$ is a semiring and $M_{\infty}$ is a tomonoid with the greatest element $+\infty$, which has been obtained from the tomonoid $M$ with no greatest element,

\item $v(xy) = v(x) + v(y)$, for all $x,y\in S$,

\item $v(x+y) \geq \min\{v(x),v(y)\}$, for all $x,y\in S$,

\item $v(1) = 0$ and $v(0) = +\infty$.

\end{enumerate}

When there is an $M$-valuation $v$ on $S$, we set $S_v = \{s\in S : v(s)\geq 0 \}$ and $P_v = \{s\in S_v : v(s) > 0 \}$. If $M = \mathbb Z$, we say that \emph{$v$ is a discrete valuation on $S$}.

\end{definition}

\begin{remark}

\label{dualdefvaluation}

Let $(S,+,\cdot)$ be a semiring and $(M,\cdot,\leq)$ a tomonoid. It is possible to give a dual definition of an $M$-valuation on $S$ as follows:

We annex an element $O$ to the tomonoid $M$ to configure the set $M_{\max}= M \cup \{O\}$ and extend ordering and monoid multiplication by the rules: $$O \leq x \text{~and~} O \cdot x = x \cdot O = O, \text{~for all~} x\in M_{\max}.$$  Imagine we can define a function $v: S \rightarrow M_{\max}$ with following properties:

\begin{itemize}

\item $v(x \cdot y) = v(x) \cdot v(y)$, for all $x,y\in S$,

\item $v(x+y) \leq \max\{v(x),v(y)\}$, for all $x,y\in S$,

\item $v(1) = 1$ and $v(0) = O$.

\end{itemize}

Then it is easy to see that $v$ defines an $M$-valuation on $S$. On the other hand, if $v$ is an $M$-valuation on $S$, by reversing the ordering, one can get a function $v: S \rightarrow M_{\infty}$ such that the three above conditions of the current remark are satisfied.

\end{remark}

\begin{remark}

\label{valuationremark}

\begin{enumerate}

\item If, in Definition \ref{valuation}, $S$ is a ring and $M$ is a totally ordered Abelian group, then our definition for valuation maps given in Definition \ref{valuation} coincides with the definition of a Bourbaki valuation on $S$ with values in $M_{\infty}$ given in \cite{Bo}. In fact, our definition for valuation maps, which is from an arbitrary semiring into a tomonoid and is based on Bourbaki's classical definition, is somehow equivalent to the definition of valuation maps defined in Definitions 2.1 and 2.3 in \cite{IKL}, that is from a semiring into a so-called commutative bipotent semiring. Also note that our definition is more general than the definition of valuation maps in \cite{J}, which is taken from \cite{IKL}, but it is only from an additively idempotent semiring into $\mathbb R \cup \{ +\infty \}$.

\item  By considering what it has been explained in Remark \ref{dualdefvaluation}, our definition for valuation maps is more general than what is defined in Definition 1.2 in \cite{T}, a paper on valuations over idempotent (i.e., characteristic one) semirings, since in the same paper, by a valuation on a semiring $S$, it is meant a surjective map $v: S \rightarrow M_{\max}$, for some totally ordered abelian group $M$, having the three conditions given in Remark \ref{dualdefvaluation} plus the following extra condition:\\

    $v(x) \leq \max\{v(x+y), v(y)\}$ for all $x,y \in S$.

\end{enumerate}

\end{remark}

Before starting to investigate the properties of valuation maps on semirings, let us give a couple of examples. Note that a semiring $S$ is called entire if $ab=0$ implies either $a=0$ or $b=0$ for all $a,b \in S$.

\begin{example}

\begin{enumerate}

\item Let $S$ be an entire semiring. One can construct a trivial $M$-valuation $v$ on $S$, by defining $v(s)=0$ for any $s\in S-\{0\}$ and $v(0) = +\infty$, where $M=\{0\}$.

\item Let $(S,+,0,\leq)$ be a tomonoid with no greatest element. Adjoin the greatest element $+\infty$ to $S$ and call the new set $S_{\infty}$. Define on $S$ the two operations $a \oplus b = \min\{a,b\}$ and $a \odot b = a+b$. One can easily check that $(S_{\infty}, \oplus,\odot)$ is a semiring and $v : S_{\infty} \rightarrow S_{\infty}$ defined by $v(s) = s$, is an $S$-valuation on $S_{\infty}$!

\item Let $\mathbb N_0$ denote the set of non-negative integers and $p\in \mathbb N_0$ be a fixed prime number. Obviously one can uniquely write any natural number $x$ in the form of $x = p^n \cdot y$, where $n\geq 0$ and $y$ is a natural number and has no factor of $p$. Now we define the map $v : \mathbb N_0 \rightarrow \mathbb N_0 \cup \{+\infty \}$ by $v(x) = v(p^n \cdot y) = n$, if $x$ is nonzero, and $v(0) = +\infty$. A simple calculation shows that $v$ is an $\mathbb N_0$-valuation on $\mathbb N_0$.

\item Let $K$ be a semifield. By a non-archimedean absolute value on $K$, we mean a function $| \cdot | : K \longrightarrow \mathbb R$, satisfying the following properties for all $x,y \in K$:

\begin{itemize}

\item $|x| \geq 0$, and $|x| = 0$ if and only if $x = 0$,
\item $|xy| = |x||y|$,
\item $|x+y| \leq \max \{|x|,|y|\}$.

\end{itemize}

It is clear that the map $v: K \longrightarrow \mathbb R \cup \{+\infty \}$, defined as $v(x) = -\ln |x|$, gives us an $\mathbb R$-valuation on $K$ (Refer to \cite{EP} and \cite{Jac}).

\end{enumerate}

\end{example}

Now let us bring the following straightforward proposition only for the sake of reference. Note that an ideal $P$ of a semiring $S$ is said to be a prime ideal of $S$, if $P \neq S$ and $ab\in P$ implies either $a\in P$ or $b\in P$ for all $a,b \in S$.

\begin{proposition}

\label{units1}

Let there exist an $M$-valuation $v$ on $S$. The following statements hold:

\begin{enumerate}

\item If $s$ is multiplicatively invertible element of $S$, then $v(s^{-1})= -v(s)$.

\item If $s^n = 1$, then $v(s) = 0$ for any $s\in S$ and natural number $n$.

\item The set $v^{-1} (+\infty)$ is a prime ideal of $S$.

\item The semiring $S$ is an entire semiring, if $v(s) = +\infty$ implies $s=0$ for all $s\in S$.

\item The set $S_v = \{s\in S : v(s)\geq 0 \}$ is a subsemiring of $S$ and $P_v = \{s\in S_v : v(s) > 0 \}$ is a prime ideal of $S_v$.

\item If $s$ is a unit of the semiring $S_v$, then $v(s) = 0$.

\end{enumerate}

\end{proposition}

The statement (5) in Proposition \ref{units1} suggests us to give the following definition:

\begin{definition}

\label{Vsemiring}

We say ``\emph{$T$ is a $V$-semiring with respect to the triple $(S,v,M)$}", if $T$ is a semiring and there exists an $M$-valuation $v$ on $S$ such that the semiring $S$ contains $T$ as a subsemiring and $T=S_v = \{s\in S : v(s)\geq 0 \} $. When there is no fear of ambiguity, we may drop the expression ``with respect to the triple $(S,v,M)$" and simply say that $T$ is a $V$-semiring.

\end{definition}

If we denote the units of a semiring $S$ by $U(S)$, then $U(S_v) \subseteq \{x\in S_v : v(x)=0 \}$ (Proposition \ref{units1}). One may ask if these two sets are equal in general. In the following, we give simple but useful examples to show that the mentioned two sets can be sometimes different:

\begin{example}

\begin{enumerate}

\item Let $S=\mathbb N_0$ be the set of all nonnegative integers and define $v: \mathbb N_0 \longrightarrow \{0, +\infty\}$ to be a map that sends all $s>0$ to $0$ and sends $0$ to $+\infty$. Then, it is easy to check that $U(S_v) = \{1\}$, while $\{s\in S_v : v(s)=0 \} = \mathbb N$.

\item Put $T=\mathbb N_0 [X][X^{-1}]$ to be the entire semiring of all Laurent polynomials on $\mathbb N_0$. Define the map $v$ from $T$ to $\mathbb Z_{\infty}$ by $v(a_m X^m + a_{m+1} X^{m+1} + \cdots + a_n X^n) = m$, where $m,n\in \mathbb Z$ and $m\leq n$ and $a_m \neq 0$ and $v(0) = +\infty$. It is easy to see that $v$ is a discrete valuation on the entire semiring $T$, $S_v = \mathbb N _0 [X]$ and $U(S_v)=\{1\}$, while $\{f\in \mathbb N _0 [X] : v(f)=0 \} = \mathbb N + X\cdot \mathbb N _0 [X]$!

\end{enumerate}

\end{example}

Now the question arises under what conditions the equality $U(S_v) = \{x\in S_v : v(x)=0 \}$ holds. This is what we are going to show in Theorem \ref{units2}:

\begin{theorem}

\label{units2}

Let $T$ be a $V$-semiring with respect to the triple $(S,v,M)$. The following statements hold:

\begin{enumerate}

\item The set $P_v = \{s\in S_v : v(s) > 0 \}$ is the only maximal ideal of the semiring $S_v$ if and only if $U(S_v) = \{x\in S_v : v(x)=0 \}$.

\item If $v$ is a surjective map, then the semiring $S$ is a semifield if and only if $M$ is an Abelian group and $U(S_v) = \{s\in S_v : v(s)=0 \}$.

\end{enumerate}

\begin{proof}
$(1)$: It is straightforward to see that $P_v = \{s\in S_v : v(s) > 0 \}$ is a prime ideal of $S_v$ (Proposition \ref{units1}). Now assume that $P_v$ is the only maximal ideal of the semiring $S_v$. We claim that any $s\in S_v$ with $v(s)=0$ is a unit. In contrary, if $s$ is a nonunit element of $S_v$, then the principal ideal $(s) \neq S_v$ is a subset of the only maximal ideal $P_v$ of $S_v$ and therefore $v(s) > 0$. On the other hand, if $U(S_v) = \{s\in S_v : v(s)=0 \}$ and $I$ is a proper ideal of $S_v$, then $I$ cannot contain a unit and this means that $v(s) > 0$ for any nonzero $s\in I$ and therefore $I \subseteq P_v$.

$(2)$: Let $M$ be an Abelian group and $U(S_v) = \{s\in S_v : v(s)=0 \}$. Let $x$ be a nonzero element of $S$ such that $v(x) = g$ for some $g\in M$. Since $v$ is surjective, there exists a nonzero $y \in S$ such that $v(y) = -g$. From this, we have $v(xy)=0$ and this means that $xy$ is unit of $S_v$ and therefore $x$ is multiplicatively invertible. On the other hand, let $S$ be a semifield. Since $v$ is surjective, $M$ is a homomorphic image of $S-\{0\}$, and the monoid $M$ is in fact a group. Also if $x\in S_v$ is a unit, then $v(x)=0$ (Proposition \ref{units1}). Now let $x\in S$ and $v(x)=0$. So there exists a $y \in S$ such that $xy=1$. This implies that $v(y)=0$ and the inverse of $x$ is in fact in $S_v$ and so $x \in U(S_v)$.
\end{proof}

\end{theorem}

\begin{remark}

\label{semifieldoffractions}

Let us recall two important points:

\begin{enumerate}

\item[(\ref{semifieldoffractions}.a)] From commutative semigroup theory, we know that any cancellative monoid $M$ can be embedded into an Abelian group $\gp(M)$, known as the group of the differences of $M$, in this sense that there is a monoid monomorphism $\iota : M \rightarrow \gp(M)$ defined by $\iota(m)=[(m,0)]_{\sim}$ such that any element $g\in \gp(M)$ can be written in the form of $g=\iota(x) - \iota(y)$ for some $x,y\in M$ \cite[p. 50]{BG}. Now let $M$ be a cancellative tomonoid. For simplification, we denote every element of $\gp(M)$ by $(x-y)$, where $x,y\in M$. The total order on $M$ induces a total order on $\gp(M)$ as follows:

    We define the relation $\leq_{gp}$ on $\gp(M)$ by $(x_1 - x_2) \leq_{gp} (y_1 - y_2)$, if $x_1 + y_2 \leq x_2 + y_1$. One can easily check that $(\gp(M), \leq_{gp})$ becomes a totally ordered Abelian group such that its order is preserved by the monomorphism $\iota$ and particularly $0 \leq x$ if and only if $0 \leq_{gp} \iota(x)$ for each $x\in M$.

\item[(\ref{semifieldoffractions}.b)] Let $S$ be a multiplicatively cancellative semiring. The semiring $S$ can be embedded into a semifield $F(S)$, known as the semifield of fractions of $S$. Note that there is a semiring monomorphism $\varepsilon : S \rightarrow F(S)$ such that any element $f\in F(S)$ can be considered in the form of $f= \varepsilon(a)\cdot \varepsilon(b)^{-1}$ with $a,b \in S$ and $b\neq 0$ \cite[p. 20]{Go1}. For simplification, we denote every element of $F(S)$ by $a/b$, where $a\in S$ and $b\in S-\{0\}$.

\end{enumerate}

\end{remark}

Now by considering Remark \ref{semifieldoffractions}, we have the following:

\begin{theorem}

\label{extension}

Let $T$ be a $V$-semiring with respect to the triple $(S,v,M)$, where $S$ is an MC semiring and $M$ is a cancellative tomonoid. Then the following statements hold:

\begin{enumerate}

\item The map $v^{\prime} : F(S) \rightarrow \gp(M)_{\infty}$ defined by $v^{\prime}(x/y) = (v(x) - v(y))$ for nonzero elements $x,y \in S$ and $v^{\prime}(0/1) = +\infty$ is a $\gp(M)$-valuation on $F(S)$.

\item The $V$-semiring $S_v$ can be considered as a subsemiring of the $V$-semiring $S_{v^{\prime}}$.

\end{enumerate}

\begin{proof}

(1): A simple calculation shows that  $v^{\prime}$ is a well-defined map from $ F(S)$ to $ \gp(M)_{\infty}$ such that $v^{\prime}((a/b)\cdot (c/d)) = v^{\prime}(a/b)+ v^{\prime}(c/d)$ and $v^{\prime}(1/1) = 0$.

Now consider that

\centerline{$v^{\prime}((a/b)+(c/d)) = v^{\prime}((ad+bc)/bd) = (v(ad+bc) - v(bd))$}
\centerline{$\geq (\min\{v(a)+v(d) , v(b)+v(c)\} - (v(b)+v(d)))$}
\centerline{$= \min\{(v(a)-v(b)) , (v(c)-v(d)) \} = \min\{v^{\prime}(a/b) , v^{\prime}(c/d) \}$.}

So we have already proved that $v^{\prime}$ is the $\gp(M)$-valuation on the semifield $F(S)$.

(2): It is easy to see that the map $\varepsilon_{\mid_{S_v}} : S_v \rightarrow \{z/1\in F(S) : z\in S_v\} $, defined by $\varepsilon_{\mid_{S_v}} (z) = z/1$, is the semiring isomorphism that we need. On the other hand, $\{z/1\in F(S) : z\in S_v\}$ is a subsemiring of $S_{v^{\prime}}$, since $v^{\prime}(z/1) = v(z)\geq 0$ for any $z\in S_v$ and this means that $z/1 \in S_{v^{\prime}}$. Therefore $S_v$ can be considered as a subsemiring of $S_{v^{\prime}}$ and this finishes the proof.
\end{proof}

\end{theorem}

\begin{remark}

The reader, who is familiar with valuation ring theory, knows that if $G$ is a totally ordered Abelian group and $\kappa$ is a field and $v : \kappa \rightarrow G_{\infty}$ is a valuation map, then $v$ has this property that $v(x) \neq v(y)$ implies $v(x+y) =  \min\{v(x),v(y)\}$ for all $x,y\in \kappa$ (\cite[Statement 1.3.4, p. 20]{EP}). Though in Remark \ref{min-property}, by giving a suitable example, we will show that, in general, this property does not hold for valuation maps on semirings, but later in Corollary \ref{subtractive2} and Theorem \ref{min-property2}, we will see that this property holds for some important families of semirings. Therefore, it is justifiable to give a name to this property.

\end{remark}

\begin{definition}

\label{min-property-def}

Let there exist an $M$-valuation $v$ on $S$. We define $v$ to have $\min$-property, whenever $v(x) \neq v(y)$, then $v(x+y) = \min\{v(x),v(y)\}$ for any $x,y \in S$.

\end{definition}

In the following remark, we give some examples related to $\min$-property for valuation maps on semirings. The first example shows that a valuation map, may not satisfy the $\min$-property. The second and third examples are good examples of valuation maps on entire semirings satisfying $\min$-property.

\begin{remark}

\label{min-property}

\begin{enumerate}

\item Let $T$ be an entire semiring. Consider the entire semiring $T[X][X^{-1}]$ of Laurent polynomials and define the map $\deg :  T[X][X^{-1}] \rightarrow \mathbb Z_{\infty}$ by $\deg(a_m X^m + a_{m+1} X^{m+1} + \cdots + a_n X^n) = n$, where $m,n\in \mathbb Z$ and $m\leq n$ and $a_n \neq 0$ and $v(0) = +\infty$. It is easy to check that the map $``\deg"$ is surjective with this property that $\deg(fg) = \deg(f) + \deg(g)$. Let us mention that a semiring $S$ is called zerosumfree, if $x+y = 0$ implies $x=y=0$ for all $x,y \in S$. Our claim is that $\deg (f+g) \geq \min \{ \deg(f), \deg(g) \}$ for $f,g \in T[X][X^{-1}]$ if and only if $T$ is a zerosumfree entire semiring and the proof is as follows:

    If $T$ is zerosumfree, then one can easily see that for all $f,g \in T[X][X^{-1}]$, we have $\deg (f+g) = \max \{ \deg(f), \deg(g) \} \geq \min \{ \deg(f), \deg(g) \}.$ And if $T$ is not a zerosumfree semiring, then there exists two nonzero elements $a,b \in T$ such that $a+b = 0$ and if one sets $f=1+aX$ and $g=bX$, then obviously $\deg(f)=\deg(g)=1$, while $\deg(f+g)=0$. From all we said we get that if $T$ is a zerosumfree entire semiring, then $``\deg"$ is a $\mathbb Z$-valuation on the entire semiring $T[X][X^{-1}]$, while the  $\min$-property does not hold for $``\deg"$.\\

    Now we give two important examples of valuation maps on entire semirings satisfying the $\min$-property:

\item Let us recall that if $M$ is a commutative monoid and $T$ is a semiring, one can define the monoid semiring $T[M]$ constructed from the monoid $M$ and the semiring $T$ similar to the standard definition of monoid rings. We write each element of $f \in T[M]$ as polynomials $f = t_1 X^{m_1}+ \cdots + t_n X^{m_n}$, where $t_1, \ldots , t_n \in T$ and $m_1, \ldots, m_n$ are distinct elements of $M$. Note that this representation of $f$ is called the canonical form of $f$ \cite[p. 68]{G2}. One can easily check that if $M$ is a tomonoid and $T$ is an entire semiring, then the map $v : T[M] \rightarrow M$ defined by $v(s_1 X^{m_1}+ \cdots + s_n X^{m_n}) = \min \{m_1, \ldots, m_n\}$, if $s_1 X^{m_1}+ \cdots + s_n X^{m_n} \neq 0$ and $v(0) = +\infty$ is an $M$-valuation on $T[M]$ satisfying the $\min$-property. This example can be interesting in another perspective because $v(f) = v(g)$ implies $v(f+g) = v(f)$, for all $f,g\in T[M]$ if and only if $T$ is a zerosumfree semiring.

\item For any given entire semiring $T$, define the entire semiring of Laurent power series $T[[X]][X^{-1}]$ to be the set of all elements of the form $\sum_{n\geq m}^{\infty} a_n X^n$, where $m\in \mathbb Z$ and $a_i \in T$. It is, then, easy to check that the map $v$, defined by $v(\sum_{n\geq m}^{\infty} a_n X^n) = m$ if $a_m \neq 0$ and $v(0) = +\infty$ on $T[[X]][X^{-1}]$, is the $\mathbb Z$-valuation on the entire semiring $T[[X]][X^{-1}]$ satisfying the $\min$-property. Note that $v(f) = v(g)$ implies $v(f+g) = v(f)$, for all $f,g\in T[[X]][X^{-1}]$ if and only if $T$ is a zerosumfree semiring.

\end{enumerate}

\end{remark}

Let $S$ be a semiring and $(M,+,0)$ be a commutative additive monoid. The monoid $M$ is said to be an $S$-semimodule if there is a function, called scalar product, $\lambda: S \times M \longrightarrow M$, defined by $\lambda (s,m)= sm$ such that the following conditions are satisfied:

\begin{enumerate}
\item $s(m+n) = sm+sn$ for all $s\in S$ and $m,n \in M$;
\item $(s+t)m = sm+tm$ and $(st)m = s(tm)$ for all $s,t\in S$ and $m\in M$;
\item $s\cdot 0=0$ for all $s\in S$ and $0 \cdot m=0$ and $1 \cdot m=m$ for all $m\in M$.
\end{enumerate}

A nonempty subset $N$ of an $S$-semimodule $M$ is said to be an $S$-subsemimodule of $M$, if $m+n \in N$ for all $m,n \in N$ and $sn \in N$ for all $s \in S$ and $n \in N$. An $S$-subsemimodule $K$ of an $S$-semimodule $N$ is said to be subtractive if $x+y \in K$ and $x\in K$ imply that $y\in K$ for any $x,y \in N$. For more on semimodules over semirings, one may refer to \cite[Chap. 14]{Go2}. The following theorem introduces some subtractive ideals of $V$-semirings. Also we will use this theorem to find some equivalent conditions for a valuation map to satisfy the $\min$-property.

\begin{theorem}

\label{subtractive1}

Let $M$ be a tomonoid and $v$ be an $M$-valuation on the semiring $T$. The following statements hold:

\begin{enumerate}

\item For any $\alpha \in M$, the sets $K_{\alpha} = \{x\in T : v(x) > \alpha \}$ and $L_{\alpha} = \{x\in T : v(x) \geq \alpha \}$ are both $S_v$-subsemimodules of $T$ and if $v$ is a surjective map and $\beta > \alpha$ is another element of $M$, then $K_{\beta} \subset L_{\beta} \subset K_{\alpha} \subset L_{\alpha}$.

\item If $\alpha \geq 0$ is an element of $M$, then $I_{\alpha} = \{x\in S_v : v(x) > \alpha \}$ and $J_{\alpha} = \{x\in S_v: v(x) \geq \alpha \}$ are both $S_v$-ideals of $S_v$ and if $v$ is a surjective map and $\beta > \alpha$ is another element of $M$, then $I_{\beta} \subset J_{\beta} \subset I_{\alpha} \subset J_{\alpha}$.

\item For any $x,y \in T$ ($\in S_v$), either $L_{v(x)} \subseteq L_{v(y)}$ ($I_{v(x)} \subseteq I_{v(y)}$) or $L_{v(y)} \subseteq L_{v(x)}$ ($I_{v(y)} \subseteq I_{v(x)}$).

\item If the valuation map $v$ has the $\min$-property, then $K_{\alpha}$ and $L_{\alpha}$ ($I_{\alpha}$ and $J_{\alpha}$) are both subtractive $S_v$-subsemimodules of $T$ (ideals of $S_v$) for any (nonnegative) $\alpha \in M$.

\end{enumerate}

\begin{proof}
The assertions (1), (2) and (3) are straightforward.

For (4), we only prove one of the claims, since the proof of the other ones is similar. Now let $x,y \in T$ be such that $x+y \in K_{\alpha}$ and $x\in K_{\alpha}$. If $v(x) = v(y)$, then from the assumption $x\in K_{\alpha}$, we get that $y\in K_{\alpha}$. If $v(x) \neq v(y)$, then from the assumptions $x+y \in K_{\alpha}$ and $v(x+y) =  \min\{v(x),v(y)\}$ ($\min$-property) and this point that $v(y) \geq  \min\{v(x),v(y)\}$, we get that again  $y \in K_{\alpha}$ and this proves that $K_{\alpha}$ is subtractive and the proof is complete.
\end{proof}

\end{theorem}

\begin{corollary}

\label{subtractive2}

Let $T$ be a $V$-semiring with respect to the triple $(K,v,M)$, where $K$ is a semifield and $v$ is a surjective map. The following statements are equivalent:

\begin{enumerate}

\item The valuation map $v$ has the $\min$-property;

\item For each nonnegative element $\alpha \in M$, the $K_v$-ideal $I_{\alpha}$ is subtractive;

\item The maximal $K_v$-ideal $P_v = \{x\in S_v : v(x) > 0 \}$ is subtractive.

\end{enumerate}

\begin{proof}

By Theorem \ref{subtractive1}, (1) implies (2) and obviously (2) implies (3). For proving $(3) \Rightarrow (1)$, we proceed as follows:

Since $K$ is a semifield, by Theorem \ref{units2}, $M$ is an Abelian group, $U(K_v) = \{x\in K_v : v(x)=0 \}$, and $P_v = \{x\in K_v : v(x) > 0 \}$ is the only maximal ideal of $K_v$. Now let $x,y \in K-\{0\}$ be such that $v(x) < v(y)$ and $g = v(y)-v(x)$. Obviously $g$ is a positive element of $M$ and since $v$ is a surjective map, there exists a nonzero $z\in K_v$ such that $v(z) = g$. This implies that $v(yx^{-1}z^{-1})=0$ and so $yx^{-1}z^{-1} \in U(K_v)$. Therefore there is a $u\in U(K_v)$ such that $y=xzu$. Now we consider that $x+y = x(1+zu)$. Obviously $zu\in P_v$ and since $P_v$ is subtractive, $1+zu$ needs to be a unit and $v(x+y) = v(x)$. This means that $v$ has the $\min$-property and the proof is complete.
\end{proof}

\end{corollary}

\begin{proposition}

\label{cyclic}

Let $T$ be a $V$-semiring with respect to the triple $(S,v,M)$. Then the following statements hold:

\begin{enumerate}

\item[(a)] For any $x\in S$, the cyclic $S_v$-subsemimodule (x) of $S_v$-semimodule $S$ is a subset of $L_{v(x)}$.

\item[(b)] For any $x\in S_v$, the principal $S_v$-ideal (x) of the semiring $S_v$ is a subset of $J_{v(x)}$.

\end{enumerate}

Moreover if $S$ is a semifield and the map $v$ is surjective, then the following statements hold:

\begin{enumerate}

\item[(c)] For any $x\in S$, the cyclic $S_v$-subsemimodule (x) of $S_v$-semimodule $T$ is equal to $L_{v(x)}$.

\item[(d)] For any $x\in S_v$, the principal $S_v$-ideal (x) of the semiring $S_v$ is equal to $J_{v(x)}$.

\end{enumerate}

\begin{proof}

The proof of (a) and (b) is straightforward. For the proof of (c), let $S$ be a semifield and $v$ surjective. By Theorem \ref{units2}, the monoid $M$ is in fact an Abelian group and $U(S_v) = \{x\in S_v : v(x)=0 \}$. If $x=0$, then there is nothing to prove. So let $x \neq 0$ and $y\in L_{v(x)}$. So $v(y) \geq v(x)$. Since $M$ is an Abelian group, we can choose $m\in M$ such that $m = v(y)-v(x)$ and since $v$ is a surjective map, there is a nonzero $z\in T$ such that $m=v(z)$. From this, we get that $v(y)-v(x)-v(z)=0$ and so $v(yx^{-1}z^{-1})=0$. This means that $yx^{-1}z^{-1} = u$ is a unit element of $S_v$ and $y=uzx$. Since $v(uz)\geq 0$, we have $y\in (x)$ and this completes the proof. The proof of (d) is similar to the proof of (c) and therefore, it is omitted.
\end{proof}

\end{proposition}

In the following example, we show that the condition of being a semifield for $S$ in the statements $(c)$ and $(d)$ in Proposition \ref{cyclic} is necessary, in this sense that if $S$ is not a semifield, the equality $(x) = J_v(x)$ may not hold:

\begin{example}

Let $\mathbb N_0$ be the set of nonnegative integers. Obviously any natural number $x > 0$ can be written uniquely in the form of $x = 5^n \cdot y$, where $n\geq 0$ and $y\in \mathbb N$ is relatively prime to $5$. We define $v_5 : \mathbb N_0 \rightarrow \mathbb N_0 \cup \{+\infty\} $ by $v_5(x) = v_5(5^n \cdot y)=n$, if $x>0$ and $v(0) = +\infty$. It is straightforward to check that $v_5$ is an $\mathbb N_0$-valuation on the semiring $\mathbb N_0$. Obviously $v(2)=0$ and $J_{v(2)} = \{ y\in \mathbb{N} : v(y) \geq 0\} = \mathbb N_0$, while $(2) = 2 \mathbb N_0$.
\end{example}

\begin{proposition}

\label{principal}

Let $T$ be a $V$-semiring with respect to the triple $(K,v,M)$, where $K$ is a semifield and $v$ is surjective. Then for the cyclic $S_v$-subsemimodules $L_1,L_2$ of the $S_v$-semimodule $K$, either $L_1 \subseteq L_2$ or $L_2 \subseteq L_1$. In particular for the principal $S_v$-ideals $J_1,J_2$ of the semiring $S_v$, either $J_1 \subseteq J_2$ or $J_2 \subseteq J_1$.

\begin{proof}
If one of the cyclic subsemimodules $L_1$ and $L_2$ is zero, then there is nothing to prove. If not, then there are two nonzero elements $x_1,x_2\in K$ such that $L_1=(x_1)$ and $L_2=(x_2)$. Since $K$ is a semifield and $v$ is surjective, by Proposition \ref{cyclic}, $L_1 = L_{v(x_1)}$ and $L_2 = L_{v(x_2)}$ and by Theorem \ref{subtractive1}, either $L_1 \subseteq L_2$ or $L_2 \subseteq L_1$. The same proof works for ideals and this is what we wished to prove.
\end{proof}

\end{proposition}

Let us recall that an integral domain $D$ is said to be a valuation ring, if every element $x$ of its field of fractions $K$ satisfies this property: $x\notin D$ implies $x^{-1} \in D$ (\cite[Chap. 4]{Mat}). This is equivalent to this statement that all ideals of $D$ are totally ordered by inclusion (\cite[Proposition 5.2]{LM}). Theorem \ref{principal} states that if $T$ is a $V$-semiring with respect to the triple $(K,v,M)$, where $K$ is a semifield and $v$ is surjective, then the principal ideals of $S_v (=T)$ are totally ordered by inclusion. The question arises if this causes all ideals of $S_v$ to be totally ordered. As we will see in Theorem \ref{valuationsemiring}, the answer to this question is affirmative. Now we pass to the next section to define valuation semirings and investigate their properties.

\section{Valuation Semirings}

The main task of this section is to define the concept of valuation semirings, which is an extension of the concept of valuation rings, and show that similar facts for valuation rings hold for valuation semirings. A classical result in valuation theory shows that $D$ is a valuation ring if and only if $D$ is a $V$-semiring with respect to the triple $(F(D),v,G_{\infty})$, where $F(D)$ is the field of fractions of the domain $D$, $G$ is the quotient group of $F(D)^*$ modulo $U(D)$ and $v: F(D) \rightarrow G_{\infty}$ is defined by $v(x) = xU$, if $x\neq 0$, and $v(0) = +\infty$ (\cite[Definition 5.12 and Proposition 5.13]{LM}). Before giving our definition for valuation semirings, we need to explain one more point: From ring theory, we know that a ring $D$ is said to be an integral domain if one of the following equivalent conditions hold:

\begin{enumerate}

\item The multiplication of the ring $D$ has cancelation property, i.e., if $ab = ac$ and $a\neq 0$, then $b=c$ for all $a,b,c \in D$,

\item The ring $D$ is entire, i.e., if $ab = 0$, then either $a=0$ or $b=0$.

\end{enumerate}

In semiring theory, these two statements are not equivalent. In fact, if $S$ is an MC semiring, then $S$ is an entire semiring but the inverse is not true. For instance, the semiring $(L = \{0,1, \ldots, n\}, \max, \min)$ for $n\geq2$ is not an MC semiring, since if $u\in L-\{0,1\}$, then $\min\{u,u\} = \min\{u,1\}$, but $u\neq 1$, while obviously it is an entire semiring.

Now imagine $T$ is a $V$-semiring with respect to the triple $(S,v,M)$. The question is: Does the assumption $T$ is entire imply that the ideals of $S_v$ are totally ordered by inclusion? The following interesting example shows that this is not the case:

\begin{example}

Let $\mathbb B = \{0,1\}$ be the Boolean semiring and $\mathbb B [X][X^{-1}]$ the entire semiring of Laurent polynomials. Define the surjective map $v: \mathbb B [X][X^{-1}] \rightarrow \mathbb Z_{\infty}$ by $v(a_m X^m + \cdots + a_n X^n) = m$, when $m \leq n$ are integer numbers and $a_m \neq 0$ and $v(0) = +\infty$. Obviously $\mathbb B[X]$ is a $V$-semiring with respect to the triple $(\mathbb B [X][X^{-1}], v, \mathbb Z)$, while even its principal ideals are not totally ordered under inclusion, since the ideals $(X)$ and $(X+1)$ of $\mathbb B[X]$ are not comparable.

\end{example}

According to this introductory note and all we have seen in Theorem \ref{extension}, Corollary \ref{subtractive2}, and Proposition \ref{cyclic}, it seems the assumptions ``$K$ is a semifield" and ``$v$ is a surjective map" are useful in defining valuation semirings. Therefore, we give the following definition:

\begin{definition}

\label{valuationsemiring-def}

We define a semiring $S$ to be a \emph{valuation semiring}, if there exists an $M$-valuation $v$ on $K$, where $K$ is a semifield containing $S$ as a subsemiring, $v$ is surjective, and $S = K_v = \{s\in K : v(s)\geq 0 \}$.

\end{definition}

\begin{remark}

\label{isomorph}

\begin{enumerate}

\item In Definition \ref{valuationsemiring-def}, according to Theorem \ref{units2}, $M$ has to be an Abelian group and $S$ needs to be an MC semiring, since it is a subsemiring of the semifield $K$.

\item In Definition \ref{valuationsemiring-def}, $K$ needs to be isomorph to the semifield of fractions of $S$. Define $\varphi : F(S) \rightarrow K$ by $\varphi(a/b) = a\cdot b^{-1}$. An easy calculation shows that $\varphi$ is a semifield monomorphism. We only prove that $\varphi$ is surjective. Let $x\in K$. If $v(x) \geq 0$, then $x\in S$ and therefore $\varphi(x/1) = x$. If $v(x) < 0$, then $v(x^{-1}) = -v(x) > 0$ and $x^{-1} \in S$ and obviously $\varphi(1/(x^{-1})) = x$.

\item Any valuation semiring is a quasi-local semiring (Theorem \ref{units2}). Note that a semiring is called to be quasi-local if it has only one maximal ideal (\cite[Example 6.1]{Go2}).

\item Let $T$ be a $V$-semiring with respect to the triple $(S,v,M)$, where $S$ is an MC semiring and $M$ is a cancellative tomonoid. Then $S_v$ can be embedded into a valuation semiring (Theorem \ref{extension}).

\end{enumerate}

\end{remark}

While the proof of the following important theorem is somehow similar to its ring version, we bring its proof here only for the completeness of the paper:

\begin{theorem}

\label{valuationsemiring}

For an MC semiring $S$, the following statements are equivalent:

\begin{enumerate}

\item $S$ is a valuation semiring,

\item For any element $x\in F(S)$, either $x\in S$ or $x^{-1}\in S$,

\item For any ideals $I,J$ of $S$, either $I \subseteq J$ or $J \subseteq I$,

\item For any elements $x,y\in S$, either $(x) \subseteq (y)$ or $(y) \subseteq (x)$.

\end{enumerate}

\begin{proof}
$(1) \Rightarrow (2)$: Let $S$ be a valuation semiring. According to Definition \ref{valuationsemiring-def}, $S = K_v = \{s\in K : v(s)\geq 0 \}$. Now take $x\in F(S)$. By considering the point (2) mentioned in Remark \ref{isomorph}, if $x\notin S$, then $v(x) < 0$ and so $v(x^{-1} > 0$ and this means that $x^{-1}\in S$.

$(2) \Rightarrow (3)$: Imagine $I \nsubseteq J$ and choose $a\in I-J$. Take $b\in J-\{0\}$. So $a/b \notin S$, because if $a/b \in S$, then $(a/b) \cdot b \in B$, which is in contradiction with $a\notin J$. So by (2), $b/a \in S$ and $(b/a) \cdot a \in I$. Thus $J \subseteq I$.

$(3) \Rightarrow (4)$: Obvious.

$(4) \Rightarrow (1)$: Let $U$ be the units of the semiring $S$. Definitely $U$ is a subgroup of $F(S)^*$, the multiplicative group of nonzero elements of fraction semifield of $S$. Define $G = F(S)^* / U$ to be the quotient group of $F(S)^*$ modulo $U$ and write the operation on G additively, i.e., $aU + bU = (ab)U$. Define a relation on $G$ by $yU \leq xU$ if for the cyclic $S$-subsemimodules $(x),(y)$ of $F(S)$, we have $(x) \subseteq (y)$. It is obvious that this is an order on $G$. On the other hand the cyclic $S$-subsemimodules of $F(S)$ are totally ordered, for the principal ideals of $S$ are totally ordered and this causes $\leq$ to be a totally order on $G$. Now we prove that $G$ with the order $\leq$ is in fact a totally ordered Abelian group. Let $yU \leq xU$. So $(x) \subseteq (y)$ and obviously $(zx) \subseteq (zy)$ and this means that $zU + yU \leq zU + xU$.

Finally define $v :  F(S) \rightarrow G_{\infty}$ by $v(x) = xU$, for $x \neq 0$ and $v(0) = +\infty$. From the definition of $G$ and $v$, it is clear that $v(xy) = (xy)U = xU + yU = v(x)+v(y)$ and so $v(1) = U$, the neutral element of the group $G$.

Also let $x\in K$. So $x\in S$ if and only if $(x) = (1)$ and this implies that $x\in S$ if and only if $U \leq xU$, which means that $x\in S$ if and only if $U \leq v(x)$. Now we use this to prove that $v(x+y) \geq \min\{v(x),v(y)\}$. If either $x=0$ or $y=0$, the inequality holds vividly. So let $x,y$ are both nonzero element of $F(S)$ and $v(y) \leq v(x)$. From this, we get that $yU \leq xU$ and according to the definition of order, $(x) \subseteq (y)$. So there exists a nonzero $s\in S$ such that $x=sy$ and therefore $x/y \in S$. Obviously this implies that $x/y + 1 \in S$ and $v(x/y + 1) \geq U$. At last we get that $v(x+y) = v(y(x/y + 1) = v(y) + v(x/y + 1) \geq v(y)$ and the proof is finished.
\end{proof}

\end{theorem}

After seeing this theorem, the reader may ask if there is a semiring, which is neither a ring nor an MC semiring, but still its ideals are totally ordered by inclusion. In the following, we show that the fuzzy semiring has this property:

\begin{example}

\label{fuzzysemiring}

Let $(\mathbb I=[0,1], \max, \min)$ be the fuzzy semiring and $A$ an ideal of $\mathbb I$. It is easy to see that $a\in A$ if and only if $[0,a] \subseteq A$ for any $a\in \mathbb I$. Therefore if $\sup (A) \in A$, then $A=[0,\sup (A)]$. Now let $\sup(A) \notin A$. Obviously $A \subseteq [0,\sup(A))$. Our claim is that $[0,\sup(A)) \subseteq A$. In order to prove the claim, we assume that $x\notin A$. Then $[x,1] \cap A = \emptyset$, because if $y\in [x,1] \cap A$, then $x\leq y$ and $x = \min\{x,y\} \in A$. On the other hand we know that $a\in A$ if and only if $[0,a] \subseteq A$ for any $a\in \mathbb I$. This means that $\bigcup_{a\in A} [0,a] \subseteq A$ and so $[x,1] \cap (\bigcup_{a\in A} [0,a]) = \emptyset$. From this, we get that $\bigcup_{a\in A}([0,a] \cap [x,1]) = \emptyset$, which implies that $[0,a] \cap [x,1] = \emptyset$ for any $a\in A$. Consequently $a < x$ for any $a\in A$ and this means that $x \notin [0, \sup(A))$. Hence, we have already proved that if $A$ is any ideal of $\mathbb I$, then either $A=[0,a]$ for some $a\in A$ or $A=[0,b)$ for some $b\in A$. Finally an easy case-by-case discussion shows that all ideals of $\mathbb I$ are totally ordered. Note that $\mathbb I$ is an entire semiring, while it is not an MC semiring.

\end{example}

 We finalize this section by generalizing another important classical result in valuation ring theory, which states that any valuation ring is integrally closed (\cite[Proposition 5.5]{LM}). Actually, we prove that if the only maximal ideal of a valuation semiring $S$ is subtractive, then $S$ is integrally closed in the sense of the following definition borrowed from the paper \cite[p. 88]{DM}. Note that if $S$ is a ring, then the following definition is equivalent to the standard definition of integrally closed rings given in Definition 4.2 in the book \cite{LM} (For more on this, also check Remark \ref{integralrem}).

\begin{definition}

\label{integral}

Let $S$ be an MC semiring and $F(S)$ its quotient semifield. The element $u\in F(S)$ is said to be integral over $S$ if there exist $a_1,a_2,\ldots, a_n$ and $b_1,b_2,\ldots,b_n$ in $S$ such that $u^n+a_1 u^{n-1}+\cdots + a_n = b_1 u^{n-1}+ \cdots + b_n$. The semiring $S$ is said to be integrally closed, if the set of elements of $F(S)$ that are integral over $S$ is the set $S$ \cite[p. 88]{DM}.

\end{definition}

\begin{proposition}

\label{integralthm}

Let $S$ be a valuation semiring such that its unique maximal ideal $J$ is subtractive. Then $S$ is integrally closed.

\begin{proof}

Clearly any element of $S$ is integral over $S$. Now let $u\in F(S)$ be integral over $S$. By Definition \ref{integral}, there exist $a_1,a_2,\ldots, a_n$ and $b_1,b_2,\ldots,b_n$ in $S$ such that $$u^n+a_1 u^{n-1}+\cdots + a_n = b_1 u^{n-1}+ \cdots + b_n. \qquad \text{($E_1)$}$$

Our claim is that $u\in S$. In contrary, suppose that $u\notin S$. By Theorem \ref{valuationsemiring}, $u^{-1}\in S$. On other hand, $u^{-1}\notin U(S)$, since if $u^{-1}\in U(S)$, then $u\in U(S)$. This means that $u^{-1} \in J$. By multiplying the both sides of the equation $E_1$ by $u^{-n}$, we get $$1+a_1 u^{-1}+\cdots + a_n u^{-n} = b_1 u^{-1}+ \cdots + b_n u^{-n}. \qquad \text{($E_2)$}$$

Obviously $ a_1 u^{-1}+\cdots + a_n u^{-n}, b_1 u^{-1}+ \cdots + b_n u^{-n} \in J$ and since $J$ is subtractive, we have $1\in J$, a contradiction. Therefore, $u\in S$ and this finishes the proof.
\end{proof}

\end{proposition}

\begin{remark}

\label{integralrem}
For more on the notions ``integral elements'', ``integral closure'', and ``integrally closed rings or ideals'' in commutative algebra, one can refer to \cite[Definition 5.1.1]{WK}), \cite[Definition 5.1.24]{WK}, \cite[p. 63]{AM}, \cite[Definition 1.1.1]{SC}, \cite[Definition 10.2.1]{BH}, and \cite[Exercise 4.14]{E}.

Also in the paper \cite{Mac}, the concept of integral closure for elements and ideals in idempotent semirings is introduced, and how it corresponds to its namesake in commutative algebra is established. Finally the term ``integrally closed ideals for semirings'' in Definition 1.4 of the paper \cite{Mac} is used, but the relationship of these notions in semirings to the notion of ``integrally closed MC-semirings'' given in page 88 of the paper \cite{DM} is unclear for the author.
\end{remark}

\section{Discrete Valuation Semirings}

 Let us recall that in Definition \ref{valuation}, we defined an $M$-valuation $v$ on $S$ to be discrete, if $M=\mathbb Z$. We start this section by defining discrete valuation semirings.

\begin{definition}

\label{discrete}

We define a semiring $S$ to be a \emph{discrete valuation semiring} (DVS for short), if $S$ is a $V$-semiring with respect to the triple $(K,v,\mathbb Z)$, where $K$ is a semifield and $v$ is surjective.

\end{definition}

Now we go further to give some examples for discrete valuation semirings that are not rings. For this, we need to bring some concepts and definitions.

Recall that if $S$ is a semiring, for $a,b \in S$, it is written $a \mid b$ and said that ``$a$ divides $b$", if $b = sa$ for some $s\in S$. This is equivalent to say that $(b) \subseteq (a)$. Also it is said that $a$ and $b$ are associates if $a=ub$ for some unit $u\in U(S)$ and if $S$ is an MC semiring, then this is equivalent to say that $(a) = (b)$. A nonzero, nonunit element $s$ of a semiring $S$ is said to be irreducible if $s = s_1 s_2$ for some $s_1, s_2 \in S$, then either $s_1$ or $s_2$ is a unit. This is equivalent to say that $(s)$ is maximal among proper principal ideals of $S$. An element $p\in S-\{1\}$ is said to be a prime element, if the principal ideal $(p)$ is a prime ideal of $S$, which is equivalent to say if $p \mid ab$, then either $p \mid a$ or $p \mid b$.

An MC semiring $S$ is called a unique factorization (or sometimes factorial) semiring if the following conditions are satisfied:

\begin{enumerate}

\item[UF1] Each irreducible element of $S$ is a prime element of $S$.

\item[UF2] Any nonzero, nonunit element of $S$ is a product of irreducible elements of $S$.

\end{enumerate}

First we construct a nice and general example of a discrete valuation semiring as follows:

\begin{example}

\label{dvs}

Let $S$ be a factorial semiring and $K = F(S)$ its semifield of fractions. Let $p\in S$ be a fixed prime element of $S$. According to the definition of factorial semirings, any nonzero $x\in S$ can be uniquely written in the form of $x = p^n \cdot x_1$, where $n\geq 0$ and $x_1$ is a nonzero element of $S$ such that it has no factor of $p$. Therefore, any nonzero $x/y \in K$, can be uniquely written in the form of $x/y = p^m \cdot x_1 / y_1$, where $m\in \mathbb Z$ and $x_1$ and $y_1$ are nonzero elements of $S$, which have no factor of $p$. Now we define a map $v_p : K \rightarrow \mathbb Z_{\infty}$ with $v_p(x/y) = v_p(p^m \cdot x_1 / y_1) = m$, when $x/y$ is a nonzero element of $K$, and $v_p(0) = +\infty$. Obviously $v_p(1) = 0$ and $v_p(x/y \cdot z/t) = v_p(x/y) + v_p(z/t)$. Now let $x/y = p^m \cdot x_1 / y_1$ and $z/t = p^n \cdot z_1 / t_1$, where $m,n\in \mathbb Z$ and $x_1, y_1, z_1$ and $t_1$ all have no factor of $p$. Without loss of generality, we can suppose that $m \leq n$. Now consider the following: $$v_p(x/y + z/t) = v_p(p^m ((x_1 t_1 + p^{n-m} y_1 z_1)/ y_1 t_1) = $$ $$m + v_p(x_1 t_1 + p^{n-m} y_1 z_1) \geq m = \min \{v_p(x/y) , v_p(z/t)\}.$$

From all we said, we get that $S_p = \{ x/y \in K : v_p(x/y) \geq 0 \} = \{ x/y \in K : p \nmid y \}$ is a discrete valuation semiring. Now the question arises if there are factorial semirings, which are not rings. The first example that may come to one's mind is $\mathbb N_0$. In Proposition \ref{factorialsemiring}, we also prove that if $D$ is a Dedekind domain, then $(\Id(D), +, \cdot)$ is a factorial semiring, where by $\Id(D)$, we mean the set of all ideals of $D$. Obviously $(\Id(D), +, \cdot)$ is not a ring.

\end{example}

\begin{lemma}

\label{idealsofdvs}

If $S$ is a discrete valuation semiring, then there exists a nonzero and nonunit element $t\in S$ such that any nonzero ideal $I$ of $S$ is of the form $I = (t^n)$ for some $n\geq 0$.

\begin{proof}
Let $S$ be a discrete valuation semiring. By definition, $S$ is a $V$-semiring with respect to a suitable triple $(K,v,\mathbb Z)$, where $K$ is a semifield and $v$ is surjective. Now let $t\in S$ such that $v(t) = 1$. Our claim is that if $I$ is a nonzero proper ideal of $S$, then there exists a natural number $n$ such that $I = (t^n)$. Put $n = \min \{v(s) : s\in I\}$. Obviously $n$ is a positive integer. Now let $s\in I$ such that $v(s) = n$. It is clear that $v(s \cdot t^{-n}) = 0$ and by Theorem \ref{units2}, $s\cdot t^{-n}$ is a unit and $(t^n) \subseteq I$. On the other hand, if $s\in I-\{0\}$, then $v(s)\geq n$ by choice of $n$ and therefore $s =u t^m$, where $u$ is a unit and $m\geq n$. This causes $s\in (t^n)$. So we have already proved that any nonzero proper ideal of $S$ is of the form $(t^n)$.
\end{proof}

\end{lemma}

Before bringing the first important theorem of this section, which characterizes discrete valuation semirings, we prove the following useful lemma. Let us recall that, similar to ring theory, we say a semiring $S$ satisfies ACCP if any ascending chain of principal ideals of $S$ is stationary.

\begin{lemma}

\label{nonunit}

Let $S$ be an MC semiring, which satisfies ACCP. Then $t\in S$ is a nonunit if and only if $ \bigcap_{n=1}^{\infty} (t^n) = (0)$.

\begin{proof}
If $t$ is a unit element of $S$, then obviously $ \bigcap_{n=1}^{\infty} (t^n) = S$. Now let $t$ be a nonunit element of $S$ and $s\in \bigcap_{n=1}^{\infty} (t^n)$. So for each natural number $n$, there exists an $s_n \in S$ such that $s = s_n t^n$. This gives us the ascending chain $(s_1) \subseteq (s_2) \subseteq \cdots \subseteq (s_n) \subseteq \cdots$, which must stop somewhere, because $S$ satisfies ACCP. Therefore there is a natural number $m$ such that $(s_m) = (s_{m+1})$. This means that there is an $r\in S$ such that $s_{m+1}=r s_m$ and so we have $s = s_{m+1} t^{m+1} = r (s_m t^m) t = rts$. Obviously this implies that $s = 0$, since if $s \neq 0$, then $t$ is a unit and the proof is finished.
\end{proof}

\end{lemma}

Let us recall that a semiring is called to be quasi-local if it has only one maximal ideal (\cite[Example 6.1]{Go2}).

\begin{proposition}

\label{units3}

Let $S$ be a semiring. Then the following statements hold:

\begin{enumerate}

\item The set of all units $U(S)$ of the semiring $S$ is equal to $S - \bigcup_{\mathfrak{m}\in \Max(S)} \mathfrak{m}$, where by $\Max(S)$, we mean the set of all maximal ideals of $S$.

\item The semiring $S$ is quasi-local if and only if $S-U(S)$ is an ideal of $S$.

\end{enumerate}

\begin{proof}
$(1):$ Let $S$ be a semiring and take $U(S)$ to be the set of all units of $S$. If $s\in U(S)$, then $s$ cannot be an element of a maximal ideal of $S$. On the other hand, if $s$ is not invertible, then the principal ideal $(s)$ of $S$ is proper and by Proposition 6.59 in \cite{Go2}, $(s)$ is contained in a maximal ideal $\mathfrak{m}$ of $S$ and therefore $s\in \mathfrak{m}$.

$(2):$ If $S$ is quasi-local and $\mathfrak{m}$ is its unique maximal ideal, then by $(1)$, $S-U(S) = \mathfrak{m}$. On the other hand, if $S-U(S)$ is an ideal of $S$, then any proper ideal of $S$ is contained in $S-U(S)$, since all elements of a proper ideal cannot be unit. This implies that $S-U(S)$ is the unique maximal ideal of $S$ and the proof is complete.
\end{proof}

\end{proposition}

\begin{theorem}

\label{ACCP}

The following statements for a semiring $S$ are equivalent:

\begin{enumerate}

\item $S$ is a discrete valuation semiring,

\item $S$ is a principal ideal MC semiring possessing a unique maximal ideal $J \neq (0)$,

\item $S$ is a unique factorization semiring with a unique (up to associates) irreducible element $t$,

\item $S$ is a quasi-local MC semiring whose unique maximal ideal $J \neq (0)$ is principal and $\bigcap_{n=1}^{\infty} J^n = (0)$,

\item $S$ is a quasi-local MC semiring, which satisfies ACCP and its unique maximal ideal $J \neq (0)$ is principal.

\item $S$ is an MC semiring and there exists a nonzero and nonunit element $t\in S$ such that any nonzero ideal $I$ of $S$ is of the form $I = (t^n)$ for some $n\geq 0$.

\end{enumerate}

\begin{proof}

First we prove $(1) \Rightarrow (4) \Rightarrow (2) \Rightarrow (3) \Rightarrow (1)$.

$(1) \Rightarrow (4)$: Let $S$ be a discrete valuation semiring. It is clear that $S$ is an MC semiring. Also by Lemma \ref{idealsofdvs}, we know that $S$ is quasi-local and its unique maximal ideal $J$ is nonzero and principal generated by an element $t\in S$ such that $v(t) =1$. On the other hand, any ideal of $S$ is principal. So by \cite[Proposition 6.16]{Go2}, $S$ is Noetherian and in particular it satisfies ACCP. Now since $t$ is nonunit (because $v(t) = 1$), by Lemma \ref{nonunit}, $\bigcap_{n=1}^{\infty} J^n = (0)$.

$(4) \Rightarrow (2)$: Let $J = (t)$ be the maximal ideal of $S$ and $I$ its nonzero proper ideal. Since $\bigcap_{n=1}^{\infty} J^n = (0)$, there exists a nonnegative number $n$ such that $I \subseteq J^n$, but $I \nsubseteq J^{n+1}$. Our claim is that $J^n \subseteq I$. Let $a\in I - J^{n+1}$. So we can suppose that $a = u \cdot t^n$ for some $u\in S$, where $u \notin J$. By Proposition \ref{units3}, $u$ is a unit of $S$. This means that $(a) = (t^n) = J^n$ for any $a\in I - J^{n+1}$. Now let $x\in J^n = (a)$. So $x= s \cdot a$. But $a\in I$, so $x\in I$. This means that all nonzero proper ideals of $S$ are in the form of $J^n = (t^n)$, where $n\geq 1$.

$(2) \Rightarrow (3)$: Let $S$ be a principal ideal MC semiring such that the nonzero ideal $J = (t)$ is its unique maximal ideal. First we prove that $t$ is an irreducible element of $S$. On contrary, let $t$ be reducible, then $t=s_1 s_2$, where both $s_1$ and $s_2$ are nonunit and this causes $(t) \subset (s_1) \subset S$, which contradicts the maximality of $J$. Now let $t^{\prime}$ be another irreducible element of $S$ and $I = (m)$ be an arbitrary ideal of $S$ that contains $(t^{\prime})$. So $t^{\prime} = s\cdot m$ for some $s\in S$ and by definition either $s$ or $m$ is a unit and this means that either $(t^{\prime}) = (m)$ or $(m) = S$. From this, we get that $(t^{\prime})$ is also a maximal ideal of $S$ and by hypothesis, $(t^{\prime}) = (t)$, which means that $S$ has a unique (up to associates) irreducible element $t$. Since $S$ is a principal ideal semiring, it is Noetherian. So $S$ satisfies ACCP. Now let $s\in S$ be nonzero and nonunit. Therefore by Lemma \ref{nonunit}, there is a natural number $n$ such that $s\in (t^n)$, while $s\notin (t^{n+1})$. Obviously this implies that $s= u t^n$, where $u$ is a unit element of $S$. So $S$ is a unique factorization semiring.

$(3) \Rightarrow (1)$: Let $S$ be a unique factorization semiring with a unique (up to associates) irreducible element $t$. So any element of $S$ can be written uniquely in the form of $ut^n$ for a unit $u$ and a nonnegative integer $n$. This means that any element of the fraction semifield $F(S)$ of $S$ can be written uniquely in the form of $ut^n$, where $u$ is a unit and $n$ is an integer number. Now define a map $v : F(S) \rightarrow \mathbb Z_{\infty}$ by $v(ut^n) = n$ for any $n\in \mathbb Z$ and $v(0) = +\infty$. It is a routine discussion that $S$ is a $V$-semiring with respect to the triple $(F(S), v, \mathbb Z)$.\\

Now we show $(4) \Leftrightarrow (5) \Leftrightarrow (6)$.

$(4) \Rightarrow (5)$: It is clear that (2) implies (5) and since (4) is equivalent to (2), it means that (4) implies (5).

$(5) \Rightarrow (4)$: Let $J = (t)$ be the nonzero unique maximal ideal of $S$. Since $t$ is nonunit, by Lemma \ref{nonunit}, $ \bigcap_{n=1}^{\infty} (t^n) = (0)$ and it means that $\bigcap_{n=1}^{\infty} J^n = (0)$.

$(5) \Rightarrow (6)$: By Lemma \ref{idealsofdvs}, $(1)$ implies (6) and since (5) is equivalent to (1), it is clear that $(5)$ implies (6).

$(6) \Rightarrow (5)$: Straightforward.
\end{proof}

\end{theorem}

Let us recall that a semiring $S$ is called an Euclidean semiring if there is a function $\delta : S-\{0\} \rightarrow \mathbb N_0 \cup \{+\infty\}$ satisfying the following condition:

For $a,b \in S$ with $b\neq 0$, there exist $q,r \in S$ satisfying $a = qb + r$ with $r = 0$ or $\delta(r) < \delta(b)$ \cite[Chap. 12, p. 136]{Go2}.

\begin{proposition}

Any discrete valuation semiring is an Euclidean semiring.

\begin{proof}
Let $S$ be a $V$-semiring with respect to the triple $(K,v,\mathbb Z)$, where $K$ is a semifield and $v$ is a surjective map. Define $\delta : S-\{0\} \rightarrow \mathbb N_0 \cup \{+\infty\}$ by $\delta(s) = v(s)$ for any nonzero $s\in S$ and suppose that $a,b\in S$ and $b\neq 0$. If $\delta(a) < \delta(b)$, then $a = 0\cdot b + a$. If $\delta(a) \geq \delta(b)$ then $q = ab^{-1} \in S$ and therefore $a = qb + 0$ and this completes the proof.
\end{proof}

\end{proposition}

Now we are ready to prove another important theorem of this section. Before that we need to recall the concept of Gaussian semirings, introduced in \cite[Definition 7]{N}: Let $S$ be a semiring. For $f \in S[X]$, the content of $f$, denoted by $c(f)$, is defined to be the $S$-ideal generated by the coefficients of $f$. A semiring $S$ is called Gaussian if $c(fg)=c(f)c(g)$ for all polynomials $f,g \in S[X]$. Also note that a semiring $S$ is called a subtractive semiring if each ideal of the semiring $S$ is subtractive.

\begin{theorem}

\label{min-property2}

Let $S$ be a discrete valuation semiring. Then the following statements are equivalent:

\begin{enumerate}

\item The unique maximal ideal of $S$ is subtractive,

\item The valuation map $v : F(S) \rightarrow \mathbb Z_{\infty}$ satisfies the $\min$-property,

\item $S$ is a subtractive semiring,

\item $S$ is a Gaussian semiring.

\end{enumerate}

\begin{proof}

$(1) \Rightarrow (2)$: This is obvious by Corollary \ref{subtractive2}.

$(2) \Rightarrow (3)$: Let the valuation map $v : F(S) \rightarrow \mathbb Z_{\infty}$ satisfy the $\min$-property. Let $I$ be a nonzero proper ideal of $S$. By Lemma \ref{idealsofdvs}, there exists a $t\in S$ with $v(t) = 1$ and $n\in \mathbb N$ such that $I = (t^n) = \{x\in S : v(x) \geq n \}$ and by Theorem \ref{subtractive1}, $I$ is subtractive. So we have already proved that each ideal of $S$ is subtractive.

$(3) \Rightarrow (4)$: Let $S$ be a subtractive semiring and $f,g\in S[X]$ be two arbitrary polynomials. By Dedekind-Mertens lemma for semirings (\cite[Theorem 3]{N}), $c(f)^{m+1} c(g) = c(f)^m c(fg)$, where $m = \deg(g)$. If $f = 0$, then there is nothing to prove. So suppose $f \neq 0$. Our claim is that $c(f)$ is a nonzero principal ideal of $S$. In fact, since in valuation semirings principal ideals are totally ordered by inclusion (Theorem \ref{valuationsemiring}), each finitely generated ideal has to be principal. But $S$ is an MC semiring, so $c(f)$ is a cancelation ideal of $S$ and $c(f)c(g) = c(fg)$ (\cite[Proposition 15]{N}).

$(4) \Rightarrow (1)$: If $S$ is a Gaussian semiring, then by \cite[Theorem 3]{N}, $S$ is a subtractive semiring. Therefore any ideal of $S$, in particular the unique maximal ideal of $S$ is subtractive and this finishes the proof.
\end{proof}

\end{theorem}

\begin{example}

Let $\mathbb Q^{\geq 0}$ be the semifield of all nonnegative rational numbers. Every nonzero element $a/b \in   \mathbb Q^{\geq 0}$ can be written uniquely in the form of $p^n \cdot a_1 / b_1$, where $n\in \mathbb Z$ and $a_1,b_1$ are both natural numbers and relatively prime to the given prime number $p$. Now we define $v_p (a/b) = p^n \cdot a_1 / b_1 = n$ (Refer to Example \ref{dvs}). It is, then, easy to see that $v_p$ is a discrete valuation on $\mathbb Q^{\geq 0}$ satisfying the $\min$-property. Therefore the semiring $\mathbb Q^{\geq 0}_p = \{a/b \in \mathbb Q^{\geq 0} : p \nmid b\}$ is a simple but good example for a subtractive discrete valuation semiring.

\end{example}

Now we bring the following proposition which gives us a general example for subtractive discrete valuation semirings.

\begin{proposition}

\label{factorialsemiring}

Let $D$ be a Dedekind domain and $\textbf{m}$ an arbitrary maximal ideal of $D$. Then the following statements hold:

\begin{enumerate}

\item The semiring $\Id(D)$ is a Gaussian factorial semiring;

\item The subsemiring $\{ I/J \in F(\Id(D)) : \textbf{m} \nmid J \}$ of the semifield of fractions $F(\Id(D))$ is a subtractive discrete valuation semiring.

\end{enumerate}

\begin{proof}

$(1)$: Let $D$ be a Dedekind domain. Then by \cite[Theorem 6.16]{LM}, every nonzero proper ideal of $D$ can be written as a product of prime ideals of $D$ in one and only one way, expect for the order of the factors. Also by \cite[Corollary 6.17]{LM}, every nonzero prime ideal of $S$ is a maximal ideal. On the other hand, by \cite[Theorem 6.19]{LM}, every nonzero ideal of $D$ is invertible. From all we said we get that $(\Id(D),+,\cdot)$ is a factorial semiring. Now we prove that $\Id(D)$ is a Gaussian semiring. Let $I,J \in \Id(D)$. Since $D$ is a Dedekind domain, any nonzero ideal of $D$ is invertible. Also since $I \subseteq I+J$, by \cite[Theorem 6.20 (8)]{LM}, there exists an ideal $K$ of $D$ such that $I=K(I+J)$. So $I$ belongs to the principal ideal of the semiring $\Id(D)$, generated by $I+J$, i.e., $I \in (I+J)$. Similarly it can be proved that $J \in (I+J)$. So we have $(I,J) = (I+J)$ and by \cite[Theorem 8]{N}, $\Id(D)$ is a Gaussian semiring.

$(2)$: By considering the Example \ref{dvs} and Theorem \ref{min-property2}, the semiring $S_{\textbf{m}} = \{ I/J \in F(\Id(D)) : \textbf{m} \nmid J \}$ is a discrete valuation semiring and so we only need to prove that the valuation map $v_{\textbf{m}} : F(\Id(D)) \rightarrow \mathbb Z_{\infty}$  has $\min$-property. So let $I/J = \textbf{m}^p I_1 / J_1$ and $K/L = \textbf{m}^q K_1 / L_1$ be arbitrary elements of $F(\Id(D))$ such that $\textbf{m} \nmid I_1, J_1, K_1, L_1$. Without loss of generality one may assume that $p < q$. Therefore we have: $v_\textbf{m}( \textbf{m}^p ((I_1 L_1 + \textbf{m}^{(q-p)} J_1 K_1) / J_1 L_1) = p +  v_\textbf{m}(I_1 L_1 + \textbf{m}^{(q-p)} J_1 K_1)$. Our claim is that $\textbf{m} \nmid I_1 L_1 + \textbf{m}^{(q-p)} J_1 K_1$. In contrary, if $\textbf{m} \mid I_1 L_1 + \textbf{m}^{(q-p)} J_1 K_1$, so $I_1 L_1 + \textbf{m}^{(q-p)} J_1 K_1$ is an element of the principal ideal $(\textbf{m})$. But we have proved in above that $\Id(D)$ is a Gaussian semiring and by \cite[Theorem 3]{N}, $\Id(D)$ is a subtractive semiring and in particular $(\textbf{m})$ is a subtractive ideal of $\Id(D)$ and this implies that $I_1 L_1 \in (\textbf{m})$, which contradicts our assumption that $\textbf{m} \nmid I_1, L_1$. So $v_{\textbf{m}} (I_1 L_1 + \textbf{m}^{(q-p)} J_1 K_1) = 0$ and $v_{\textbf{m}} (I/J + K/L) = \min \{v_{\textbf{m}} (I/J), v_{\textbf{m}} (K/L) \}$ and the proof is complete.
\end{proof}

\end{proposition}

\begin{remark}

\label{idealdedekindnatural}

Let $\Lambda$ be a nonempty set. The set of all functions $f: \Lambda \rightarrow \mathbb N_0$ with finite support is denoted by $\mathbb N_0 ^{(\Lambda)}(= \bigoplus_{\lambda\in \Lambda} \mathbb N_0)$. We put $S_{\Lambda} = \mathbb N_0 ^{(\Lambda)} \cup \{+\infty\}$ and define the operations ``$\min$'' and ``$+$'' on $\mathbb N_0 ^{(\Lambda)}$ componentwise and extend these operations on $S_{\Lambda}$ as follows:

\begin{itemize}

\item $\min \{f, +\infty\} = \min \{+\infty, f\} = f$, for all $f\in S_{\Lambda}$,

\item $f + (+\infty) = (+\infty) + f = +\infty$, for all $f\in S_{\Lambda}$.

\end{itemize}

It is easy to verify that $(S_{\Lambda}, \min,+)$ is a semiring such that $+\infty$ is its zero and the constant zero function, i.e., $f(\lambda)=0$ for all $\lambda \in \Lambda$, its $1_{S_{\Lambda}}$. Our claim is that if $D$ is a Dedekind domain, then $(\Id(D),+,\cdot)$ and $(S_{\Lambda}, \min,+)$ are isomorphic as semirings for $\Lambda = \Max(D)$ and the routine proof is based on an interesting result mentioned in Proposition 6.8 in \cite{Kar}, which states that if $I = \mathfrak{m}_1 ^ {\alpha_1} \cdots \mathfrak{m}_n ^ {\alpha_n}$ and $J = \mathfrak{m}_1 ^ {\beta_1} \cdots \mathfrak{m}_n ^ {\beta_n}$, then $I+J = \mathfrak{m}_1 ^ {\min\{\alpha_1,\beta_1\}} \cdots \mathfrak{m}_n ^ {\min\{\alpha_n,\beta_n\}}$, where $\alpha_i,\beta_i$ are non-negative integers, $\{\mathfrak{m}_i\}$ are distinct maximal ideals of $D$, and $\mathfrak{m}_i ^ 0$ is taken to be the ideal $D$.

The question may arise if Proposition \ref{idealdedekindnatural} can be generalized for the semiring $(S_{\Lambda}, \min,+)$, where $S_{\Lambda} = \mathbb N_0 ^{(\Lambda)} \cup \{+\infty\}$ and $\Lambda$ is an arbitrary nonempty set. In fact, a wonderful result by Luther Claborn, mentioned and proved in Theorem 15.18 in \cite{F}, states that if $\Lambda$ is an arbitrary nonempty set, then there is a Dedekind domain $D$ with a bijection $\sigma: \Max(D) \rightarrow \Lambda$.

\end{remark}

\section*{Acknowledgements}

The author was partly supported by Department of Engineering Science at University of Tehran and wishes to thank the anonymous referee for an extensive list of comments and suggestions, which greatly helped to improve the paper. This paper is dedicated to the memory of Prof. Dr. Manfred Kudlek, since he was the professor who introduced the algebraic structure of semiring to me during my stay in Tarragona, Spain.

\bibliographystyle{plain}

\begin{thebibliography}{15.}

\bibitem{AM} M. F. Atiyah and I. G. Macdonald, {\em Introduction to Commutative Algebra}, Addison-Wesley, Reading, Massachusetts, 1969.

\bibitem{Bo} N. Bourbaki, {\em Commutative Algebra. Chapters 1--7, Elements of Mathematics}, Springer, 1998.

\bibitem{BG} W. Bruns and J. Gubeladze, {\em Polytopes, Rings, and $K$-Theory}, Springer, Dordrecht, 2009.

\bibitem{BH} W. Bruns and J. Herzog, {\em Cohen-Macaulay Rings}, 2nd edn., Cambridge studies in advanced mathematics: Vol. 39, Cambridge Univ. Press, Cambridge, 1998.

\bibitem{DM} D. J. Dulin and J. R. Mosher, {\em The Dedekind property for semirings}, J. of Aust. Math. Soc., {\bf 14} (1972), 82--90.

\bibitem{E} D. Eisenbud, {\em Commutative Algebra with a View Toward Algebraic Geometry}, Springer-Verlag, New York, 1995.

\bibitem{EP} A. J. Engler and A. Prestel, {\em Valued Fields}, Springer, Berlin, 2005.

\bibitem{F} R. M. Fossum, {\em The Divisor Class Group of a Krull Domain}, Ergebnisse der Mathematik und ihrer Grenzgebiete, Band {\bf 74}, Springer-Verlag, Berlin, 1973.

\bibitem{G2} R. Gilmer, {\em Commutative Semigroup Rings}, The University of Chicago Press, 1984.

\bibitem{Go1} J. S. Golan, {\em Power Algebras over Semirings}, Springer, 1999.

\bibitem{Go2} J. S. Golan, {\em Semirings and Their Applications}, Kluwer Academic Publishers, Dordrecht, 1999.

\bibitem{HW} U. Hebisch and H. J. Weinert, {\em Semirings - Algebraic Theory and Applications in Computer Science}, World Scientific, Singapore, 1998.

\bibitem{IKL} Z. Izhakian, M. Knebusch, and L. Rowen, {\em Supertropical semirings and supervaluations}, J. Pure Appl. Algebra {\bf 215}, No. 10 (2011), 2431--2463.

\bibitem{Jac} N. Jacobson, {\em Basic Algebra II}, 2nd ed., W. H. Freeman and Company, New York, 1989.

\bibitem{J} J. Jun, {\em Valuations of Semirings}, arXiv preprint arXiv:1503.01392 (2015).

\bibitem{Ka} I. Kaplansky, {\em Commutative Rings}, Allyn and Bacon, Boston, 1970.

\bibitem{Kar} G. Karpilovsky, {\em Topics in Field Theory}, North-Holland Mathematics Studies Vol. {\bf 155}, Amsterdam: North-Holland Publishing Company, 1989.

\bibitem{Kr} W. Krull, {\em Beitr\"{a}ge zur Arithmetik kommutativer Integrit�tsbereiche. VI. Der allgemeine Diskrimminantensatz. Unverzweigte Ringerweiterungen}, Math. Z., Vol. {\bf 45} (1939), 1--19.

\bibitem{LM} M. D. Larsen and P. J. McCarthy, {\em Multiplicative Theory of Ideals}, Academic Press, New York, 1971.

\bibitem{Mac} A. W. Macpherson, {\em On the difference between tropical functions and real-valued functions}, Preprint, arXiv:1507.00545, 2015.

\bibitem{Mat} H. Matsumura, {\em Commutative Ring Theory}, 2nd ed., Cambridge Studies in Advanced Mathematics, {\bf 8}, Cambridge University Press, Cambridge, 1989.

\bibitem{N} P. Nasehpour, {\em On the content of polynomials over semirings and its applications}, J. Algebra Appl., {\bf 15}, No. 5 (2016), 1650088 (32 pages).

\bibitem{P} G. Pilz, {\em Near-Rings, The Theory and Its Aplications}, Revised edition, North-Holland Publishing Company, Amsterdam-New York-Oxford, 1983.

\bibitem{SC} I. Swanson and C. Huneke, {\em Integral Closure of Ideals, Rings, and Modules}, London Mathematical Society Lecture Note Series {\bf 336}, Cambridge University Press, Cambridge, 2006.

\bibitem{T} J. Tolliver, {\em Extension of valuations in characteristic one}, Preprint,  arXiv:1605.06425, 2016.

\bibitem{WK} F. Wang and H. Kim, {\em Foundations of Commutative Rings and Their Modules}, Algebra and Applications: Vol. {\bf 22}, Springer Nature, Singapore, 2016.

\end{thebibliography}

\end{document}